\documentclass{article}
\usepackage[cp1251]{inputenc}
\usepackage[english]{babel}
\usepackage{latexsym,amsfonts,amsthm}
\newtheorem{lemma}{ Lemma}
\newtheorem{theorem}{ Theorem}
\newtheorem{definition}{ Definition}
\newtheorem{remark}{ Remark}

\def\Cl{{\mbox{Cl}}}
\def\Int{{\mbox{Int\ }}}

\input{psfig.sty}
\begin{document}
\large

\begin{center}{\sc Roots of 3-manifolds and cobordisms}
\end{center}

C. Hog-Angeloni\footnote{Partially supported by the INTAS  Project
"CalcoMet-GT" 03-51-3662}, S. Matveev\footnote{Partially supported
by the INTAS  Project "CalcoMet-GT" 03-51-3662  and the RFBR grant
05-01-0293-a}

\vspace{0.5cm}

\section{Introduction}

Given a set  of {\em simplifying moves} on 3-manifolds, we apply
them to a given 3-manifold $M$ as long as possible. What we get is
a {\em root} of $M$.  For us, it makes sense to consider three types of moves:
compressions along 2-spheres, proper discs and proper  annuli
having boundary circles in different components of $\partial M$.   Our  main result is
that for the above moves the root  of any 3-manifold exists and is unique. The
same result remains true if instead of manifolds we apply the moves to
 3-cobordisms of the type $(M,\partial_-M,\partial_+M)$. The only difference between moves on
 manifolds and moves on cobordisms
 is that   one boundary circle of   every annulus
participating in a  compression of a cobordism must lie in $\partial_-M$
while the other in $ \partial_+M$. We can also restrict ourselves  to
considering compressions along only spheres or   only spheres and
discs. The existence and uniqueness in the first case is well-known
 and essentially
comprise the content of the Milnor theorem on unique decomposition
of a 3-manifold into a connected sum. For the second case our result is close
to theorems about  
{\em characteristic  compression bodies}  and  about {\em cores}
 of irreducible manifolds,
 presented by F. Bonahon~\cite{Bo} and S. Matveev~\cite{Ma}, respectively.

We use Kneser 
existence \cite{Kn}, but perhaps the proof of the uniqueness  part
is easier with the method we are developing.

We point out that considering roots of cobordisms was motivated by
the paper~\cite{Ga} of R. Gadgil, which is interesting although
the proof of his main theorem seems to be incomplete.

We thank C.~Gordon, W.~Metzler, E.~Pervova, C.~Petronio, and
S. Zentner for  helpful discussions.  

The final version of the paper has been written by the second
author during his stay at MPIM Bonn. He thanks the institute for
hospitality and financial support.

\section{Moves, roots, and complexity }

We introduce several moves on 3-manifolds. In this paper all 3-manifolds are
assumed to be orientable. 
\begin{enumerate}

\item {\em $S$-move (compression along a 2-sphere)}.   Let $S$ be a  2-sphere in a
3-manifold $M$. Then we cut $M$ along $S$ and fill the two spheres
arising in this way with two 3-balls.

\item {\em $D$-move (compression along a disc)}.   Let $D$ be a proper   disc in $M$. Then we cut $M$ along $D$.

 \item {\em $A$-move (compression along an annulus)}.  Let $A$ be an   annulus
  in  $M$ such that its boundary
circles lie in different components of $\partial M$. Then we cut $M$
along $A$ and attach two plates $D_1^2\times I, D_2^2\times I$ by
identifying their base annuli $\partial D_1^2\times I,
\partial D_2^2\times I$  with the two copies of $A$, which appear
under cutting.
\end{enumerate}

\begin{definition} Let $F$ be a proper surface in a 3-manifold $M$
such that $F$ is either a sphere, or a disc, or an annulus. Then
$F$ is called {\em essential}, if one of the following holds:
\begin{enumerate}
\item $F$ is a sphere not bounding a ball;
\item $F$ is a disc such that the circle $\partial D$ is nontrivial in
$\partial M$;
\item $F$ is an incompressible annulus having boundary circles in different
components of $\partial M$.
\end{enumerate}
If $F$ is essential, then the corresponding $F$-move (i.e. the compression of
$M$ along $F$) is also
called essential.
\end{definition}

\begin{remark} \label{rm:incompressible} Later on under   essential surface
we will understand either an essential sphere, or an essential disc, or an
essential annulus. The condition that the
boundary circles of any essential annulus $A$  must lie in
different components of $\partial M$ guarantees us that $A$ is boundary
incompressible.
\end{remark}

\subsection{Roots and  complexity}\label{roots}

\begin{definition} Let  $M$  be a 3-manifold. Then
a 3-manifold  $N$   is called a {\em
root} of   $M$, if
\begin{enumerate}
\item $N$ can be obtained from $M$ by essential compressions along spheres,
discs, and annuli.
 \item $N$
admits no   further essential compressions.
 \end{enumerate}
\end{definition}

\begin{theorem} \label{th:roots} For any compact 3-manifold   $M$   the
 root of $M$ exists and is unique up to homeomorphism and
removing disjoint 3-spheres and balls.
\end{theorem}

We postpone the proof of the theorem to Section~\ref{proof}. Note
that  the condition on boundary circles of compression annuli to lie in
different components of $\partial M$ is
essential. Below we present an example of a 3-manifold $M$ with
two incompressible boundary incompressible annuli $A,B\subset M$
such that $\partial M$ is connected and  compressions of $M$ along $A$
and along $B$ lead us to two different 3-manifolds admitting no
further essential compression, i.e. to two different  ``roots''.

{\bf Example.} {\em Let $Q$ be the complement space of the figure
eight knot. We assume that the torus $\partial Q$ is equipped with a
coordinate system such that the slope of the meridian is (1,0).
Choose two pairs $(p,q)$, $(m,n)$ of coprime integers such that
$|q|,|n|\geq 2$ and $|p|\neq |m|$. Let $a$ and $b$ be corresponding curves
in $\partial Q$. Then the manifolds $Q_{p,q}$ and
$Q_{m,n}$ obtained by Dehn filling of $Q$ are not homeomorphic. 
By [Th], they are hyperbolic.
 
Consider the  thick  torus $X=S^1\times S^1 \times I$ and locate
its  exterior meridian  $\mu=S^1\times \{\ast\}\times \{ 1\}$ and
 interior  longitude  $\lambda= \{\ast\}\times S^1\times \{ 0\}$.
Then we  attach to $X$   two copies $Q',Q''$ of $Q$ as follows.
The first copy $Q'$ is attached to $X$ by identifying an annular
regular neighborhood $N(a)$ of $a$ in $\partial Q$ with an annular
regular neighborhood $N(\mu)$ of  $\mu$ in $\partial X$. The second copy
$Q''$  is attached by identifying $N(b)$ with $N(\lambda)$. Denote by
$M$ the resulting manifold $Q'\cup X\cup Q''$.

Since $Q$ is hyperbolic,  $M$ contains only two incompressible
boundary incompressible annuli $A$ and $B$, where $A$ is the
common image of $N(a)$ and $N(\mu)$, and  $B$ is the common image of
$N(b)$ and $N(\lambda)$.  It is easy to see that   compression of $M$
along $A$ gives us a disjoint union of a punctured $Q'_{p,q}$ and
a punctured $Q''$ while the compression along $B$ leaves us with a
punctured $Q'$ and a punctured $Q''_{m,n}$. After filling the
punctures (by compressions along spheres surrounding them), we get
two different  manifolds, homeomorphic to $Q _{p,q}\cup Q$ and
$Q_{m,n}\cup Q$. Since their connected components (i.e. $Q_{p,q},
Q_{m,n}, Q$) are hyperbolic, they are irreducible, boundary
irreducible and contain no essential annuli. Hence   $Q
_{p,q}\cup Q$ and $Q_{m,n}\cup Q$
 are different roots of $M$.}

  Let $M$ be a compact 3-manifold.
Let us apply to it  essential $S$-moves as long as possible. It
follows from Kneser finiteness theorem~\cite{Kn} that
 the number of possible  moves is bounded by a
constant depending on $M$ only. Denote by $s(M)$ the maximum of
these numbers taken over all chains of essential $S$-moves.

The following notion will be the main inductive parameter in our
proofs.

\begin{definition}\label{compl} Let $M$ be a 3-manifold. Then the
 complexity ${\bf c}(M)$ of $M$ is the pair
 $(g^{(2)}(\partial M),s(M))$,  where $g^{(2)}(M)=\sum_i g^2(F_i)$,   
 $g(F_i)$ is the genus
 of a component $F_i\subset \partial M$,
   and the sum
is taken over all    components of $\partial M$. The pairs are considered in   lexicographical order.
\end{definition}

The use of complexity as an inductive parameter is justified
by the following fact.

\begin{lemma}\label{lm:decrease}
Each    essential $S,D,A $-move strictly decreases
  ${\bf c}(M)$.
\end{lemma}
\begin{proof} If an essential $D$-move cuts $\partial M$ along a nonseparating curve
$\ell$ on some component $F $ of $\partial M$, then it strictly decreases $g(F)$ and hence
$c(M)$. If the move turns $F$ into two components $F', F''$, then 
$g(F)=g(F')+g(F'')$ and, since $\ell$ is nontrivial and thus $g(F')$,
$g(F'')\neq 0$, we have $g^2(F)>g^2(F')+g^2(F'')$. This implies that $c(M)$ is
decreased again. 
  The case of the $A$-move is similar.
 
 As   follows from the definition of $s(M)$,   each
 essential $S$-move strictly decreases $s(M)$. The boundary of $M$ remains
the same, hence so  does $g^{(2)}(M)$.
\end{proof}
\begin{remark} \label{rem:A-increase} It is easy to show that inessential
 $S$- and $D$-moves preserve
the complexity. However, an inessential $A$-move can increase it, but
only at the expense of $s(M)$  ($g^{(2)}(M)$ cannot increase). For
example, if an annulus $A$ cuts off a $D^2\times I$  from $M$,
then the corresponding move results in the appearing of an additional
component of the type $S^2\times I$.
\end{remark}

\subsection{Equivalence of essential surfaces}
Throughout this section, surface means sphere or disc or annulus.
\begin{definition} Let $M$ be a 3-manifold
and   $F, G$  be two
essential surfaces in $M$.  Then $F,G$ are  {\em equivalent} (we write $F\sim
G$) if
there exists a finite sequence of essential surfaces $ X_1, X_2,
\dots , X_n$ such that the following holds:

\begin{enumerate}
\item $F=X_1$ and $X_n=G$;
\item For each $i,1\leq i<n,$ the surfaces $X_i$ and $X_{i+1}$ are
disjoint.
\end{enumerate}
\end{definition}

\begin{lemma}\label{lm:all equiv}
Let $M$ be a 3-manifold   not homeomorphic to $S^1\times S^1\times I$.
Then any two
essential surfaces in $M$ are equivalent.
\end{lemma}

\begin{proof} Let $F,G$ be two essential surfaces in $M$ in general position.
Then the number of curves (circles and arcs) in the intersection
of $F$ and $G$ will be denoted by $\#(F\cap G)$. Arguing by
induction, we may assume that any two essential surfaces $F,G$  with
$\#(F\cap G)<n$  are equivalent. The base of the induction is
evident: if $\#(F\cap G)=0$, then $F\sim G$ by definition. Let
$F,G$ be two essential surfaces such that $\#(F\cap G)=n$.

{\em Case 1.  } Suppose that  $F\cap G$ contains a circle  $s$ which is trivial
in $F$. By an innermost circle argument we may assume that $s$ bounds a disc
$D\subset F$  such that $D\cap G=s$.   Compressing $G$ along $D$, we get
a two-component surface $G'$ such that one component is a sphere, the other
  is homeomorphic to $G$,
and $\#(F\cap G')=n-1$. Since $G$ is an interior connected sum of
 the components of $G'$,  at least one of them (denote it by $X$) is
essential and thus   $F\sim X$ by the inductive assumption.
On the other hand, $X$ can be shifted away from  $G$ by a small
isotopy.  It follows that $X\sim G$ and thus $F\sim G$.

 {\em Case 2.  } Suppose that  $F\cap G$ does not contain  trivial circles,
 but contains an arc $a$ which is trivial in $F$. By an outermost arc argument
 we may assume that $a$ cuts off a disc
$D\subset F$  from $F$ such that $D\cap G=a$.
Compressing $G$ along $D$, we get a
two-component surfaces $G'$ such that one component is a proper disc, the other
  is homeomorphic to $G$,
and $\#(F\cap G')=n-1$. Since $G$ is an interior boundary
connected sum of the components of $G'$,  at least one of them (denote
it by $X$) is essential and thus equivalent to $F$ by the
inductive assumption. On the other hand, $X$ can be shifted away from
$G$ by a small isotopy.  It follows that $G\sim X$ and thus
$F\sim G$.

 {\em Case 3.} Suppose that $F$ and $G$ are annuli such that $F\cap G$ consists of circles parallel to the
 core circles of $F$ and $G$. Then one can find two different components $A, B$ of $\partial M$ such that
 a circle of $\partial F$ is in $A$ and a circle of $\partial G$ is in $B$. Denote by $s$
  the first circle of $F\cap G$ we meet at our radial way along $F$ from the
  circle $\partial F\cap A$ to the other boundary circle of $F$.
  Let $F'$ be the subannulus
  of $F$ bounded by   $\partial F\cap A$ and $s$, and
     $G'$   the subannulus of $G$ bounded by     $s$ and $\partial G\cap B$.
      Then the
 annulus $F'\cup G'$ is essential and is isotopic to an 
 annulus $X$ such that $\#(X\cap F)<n$ and
 $\#(X\cap G)=0$, see Fig.~\ref{newannulus} (to get a real picture, multiply by
 $S^1$). It follows that $F\sim G$.

  \begin{figure}
\centerline{\psfig{figure=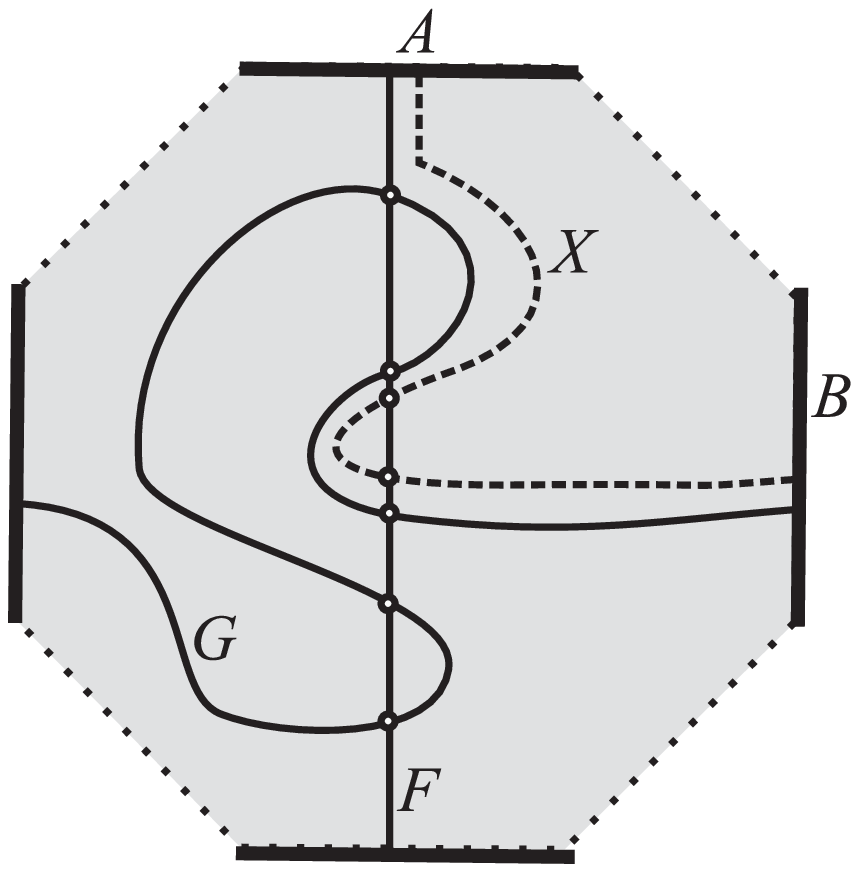,height=5.6cm}}
  \caption{ }
  \label{newannulus}
\end{figure}

  {\em Case 4.} Let  $F$ and $G$ be annuli such that $F\cap G$ consists
  of more than one  radial segments, each having endpoints in
  different  components of $\partial F$ and
  different  components of $\partial G$.

  {\em Case 4.1.}
  Suppose that there are two neighboring segments
   $s_1,s_2  \subset F\cap G\subset F$ such that
    $G$ crosses $F$ at  $s_1,s_2$  in opposite directions. Denote
    by $D$ the
quadrilateral part  of $F$
     between them. Then we cut
        $G$ along $s_1,s_2$ and attach to it   two parallel copies of $D$
         lying on different
        sides of $F$. We get a new surface $G'$
          consisting of two disjoint annuli, at least
 one of which (denote it by $X$) is essential, see Fig.~\ref{lick}
 to the left. Since $\#(X\cap
 F)=n-2$ and, after a small isotopy of $X$, $\#(X\cap
 G)=\emptyset$, we get   $F\sim X\sim  G$.

   \begin{figure}
\centerline{\psfig{figure=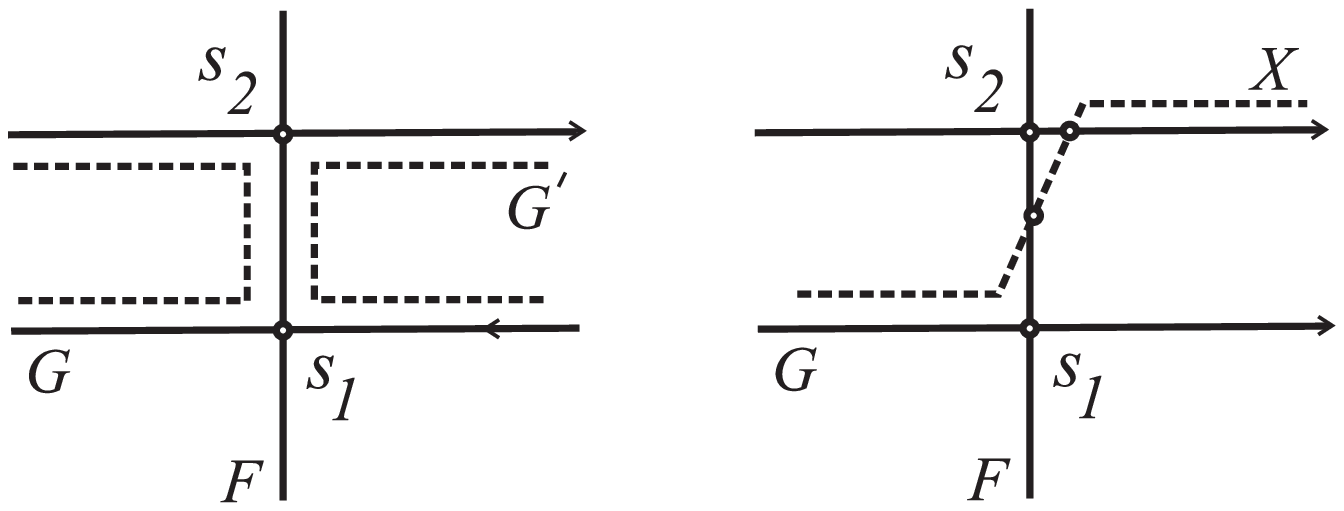,height=3.6cm}}
  \caption{ }
  \label{lick}
\end{figure}

{\em Case 4.2.} Suppose that at all segments $G$ crosses $F$ in
the same direction (say, from the left to the right). Let $s_1,
s_2$ be two neighboring segments spanned by a quadrilateral part
$D\subset F$ between them. Then $s_1,s_2$ decompose $G$ into two
strips $L_1,L_2$ such that $L_1$ approaches  $s_1$ from the left
side   of $F$ and   $s_2$ from the  right   side. Then the
annulus $L_1\cup D$ is isotopic to an annulus $X$ such that
$\#(X\cap F)\leq n-1$ and $\#(X\cap G)=1$, see Fig.~\ref{lick}  to
the right. Since $X$ crosses $F$ one  or more times in the same
direction, it is essential.
 Therefore, $F\sim X\sim  G$.

{\em Case 5.} This is the last logical possibility. Suppose that
$F$ and $G$ are annuli such that $F\cap G$ consists of one radial
segment. Denote by $G'$ the relative boundary $\partial_{rel}(N)=\Cl(N\cap
\Int M)$ of a regular neighborhood $N$ of $F\cup G$ in $M$. Then
$G'$ is an annulus having boundary circles in different components
of $\partial M$.

  {\em Case 5.1.} If $G'$ is incompressible, then we put $X=G'$.

{\em Case 5.2.}   If $G'$ admits a compressing disc $D$,   then
the relative boundary of a regular neighborhood $N$  of $G'\cup D$
consists of a parallel copy of $G'$ and two proper discs
$D_1,D_2$. If at least one of these discs (say, $D_1$) is
essential, then we put $X=D_1$.

 {\em Case 5.3.} Suppose that   discs $D_1, D_2$ are not
 essential. Then the circles $\partial D_1, \partial D_2$ bound discs
 $D_1',D_2'$ contained in the corresponding components of $\partial M$.  We
 claim that at least one of the spheres $S_1=D_1\cup D_1',S_2=D_2\cup D_2'
 $ (denote it by $X$)  must be essential. Indeed, if both bound balls, then $M$ is
 homeomorphic to $S^1\times S^1\times I$, contrary to  our  assumption.

 In all three cases 5.1-5.3 $X$ is disjoint to $F$ as well as
 to $G$. Therefore,  $F\sim X\sim  G$.
 \end{proof}

 \section{Proof of the main theorem}
 \label{proof}

Let $F$ be a sphere, a disc or an annulus in a 3-manifold $M$.
It is convenient to denote by $C_F(M)$ the result of the   $F$-move, i.e.
 the manifold obtained by     compressing  $M$ along $F$.

\begin{figure}
 \centerline{\psfig{figure=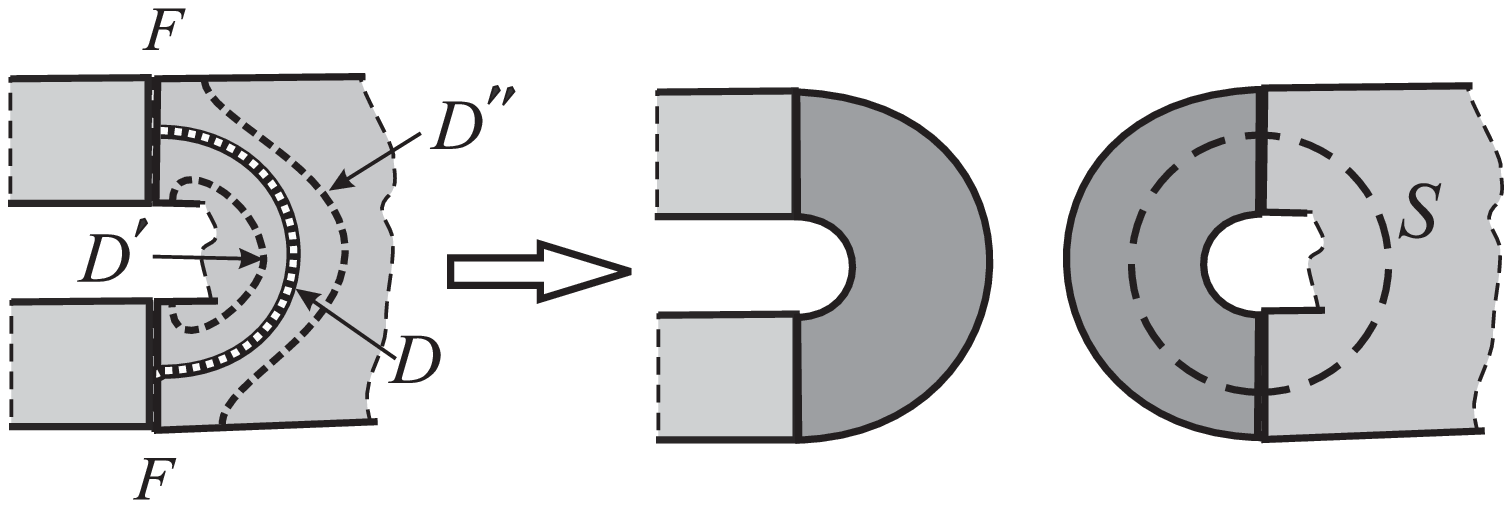,height=3.3cm}}
   \caption{ }
   \label{ddas}
 \end{figure}

 \begin{lemma}\label{commonroot} If $F$ is a sphere or a disc or an essential
 annulus, then any
 root of $M_F=C_F(M)$ is a root of $M$. If $F$ is an inessential annulus, then
   $M_F$ and $M$ have at least one common root.
\end{lemma}
 \begin{proof} It is convenient to decompose the proof into four steps.
 \begin{enumerate}
 \item[(1)] If $F$ is essential, then   any
 root of $M_F$ is a root of $M$ by definition of the root.

 \item[(2)] If $F$ is an inessential sphere,
 then $M_F$ is a union of $M$ and a disjoint 3-sphere. Therefore,
all roots of $M$ and $M_F$ are the same.

 \item[(3)] Let  $F$ be an inessential disc. Then its boundary circle bounds a
 disc $D\subset \partial M$. Choose  a 2-sphere $S$ inside $ M$ which is
  parallel to the 2-sphere $F\cup D$. Then the manifold
  $M_F$ is obtained from the manifold $M_S=C_S(M)$ by puncturing (cutting off a
  ball $V\subset M_S$). We claim that any root $R$  of
   $M_F  = M_S\setminus \Int V\subset M_S$, which can be obtained from $M_F$ by
  successive compressions along essential subsurfaces,
   is a root of $M_S$. Indeed,
  we simply compress $M_S$ along the same surfaces and get either $R$ (if
    one of those subsurfaces is a sphere surrounding $V$) or a punctured $R$
    (if the puncture survives all compressions). One more compression along a
    sphere surrounding the puncture is sufficient to convert the punctured $R$
    to $R$ (modulo disjoint 3-spheres and balls, which are irrelevant).
     It follows from (1), (2), 
     and (3)   that any root of
    $M_F$ is a root of $M$.
 \item[(4)] Let $F$ be an inessential (i.e. compressible) annulus and $D $ a
 compressing disc for $F$ such that $D\cap F=\partial D$ is a core circle of $F$.
 Denote by $N$ a regular neighborhood of $F\cup D$ in $M$. Then the
  relative boundary $\partial_{rel}N=\Cl(\partial N\cap \Int M)$ consists of a parallel
 copy of $F$ and two proper discs $D', D''$. Denote by $S$ a 2-sphere in
 $C_F(M)$ composed from a copy of $D$ and a core disc of one of the attached
 plates, see Fig.~\ref{ddas}. Then the manifolds
 $C_{D''}(C_{D'}(M))$ and $ C_S(C_F(M))$ are
 homeomorphic. Applying (1) - (3), we conclude that any root of the manifold
 $C_{D''}(C_{D'}(M))= C_S(C_F(M))$  is a root of both $M$ and $C_F(M)$.
\end{enumerate}  \end{proof}

 Proof of Theorem 1. {\sc Existence.} Let us apply to $M$ essential $S,D,A$-moves
 in arbitrary order as long as possible. By Lemma~\ref{lm:decrease},
 each move strictly decreases  the complexity. Since every set of pairs
 of nonnegative integers has a minimal pair, the process stops and we get a
 root.

 {\sc Uniqueness.} Assume the converse: suppose that there is a
 3-manifold   having two different roots. Among all
 such manifolds we choose a manifold $M$ having minimal
 complexity. 
 Then there exist two sequences of essential moves  
producing two different roots. Denote by $C_F$ and $C_G$ the first moves of
 the sequences, where $F,G$ are essential surfaces in $M$. 
 By Lemma~\ref{lm:all equiv}, there are essential surfaces
$ X_1, X_2, \dots , X_n$ such that   $F=X_1$, $X_n=G$, and that the 
 surfaces  $X_i$ and $X_{i+1}$ are
disjoint for all  $i,1\leq i<n$.  We may begin the construction of a
root starting with the compression along any of them. Evidently, for
at    least two neighboring surfaces $X_k, X_{k+1}$ the roots thus
obtained are different. For convenience, we rename $X_k,X_{k+1}$
by $F,G$ thus getting two disjoint surfaces such that $C_F(M)$ and $C_G(M)$
have different roots.  
  Then  $F$ is a subsurface
 of $M$ and of $M_G=C_G(M)$ while $G$ is a subsurface of $M$ and of $C_F(M)$.
 Denote by $N$ the manifold, obtained from $M$ by compressions along
 both surfaces $F,G$. Of course, it coincides with    $C_G(C_F(M))$  and
 $ C_F(C_G(M))$.
 
 We claim that the complexity of $N$ is strictly less than the one of $M$.
 Indeed, if $F$ is either a sphere or a disc, then  $c(N)\leq c(M_G)$
  (since compression along a sphere or a disc does not
 increase complexity) while $c(M_G)<c(M)$ by Lemma~\ref{lm:decrease}.
 Suppose that $F$ is an annulus. Then     $g^{(2)}(\partial N)$ is no greater than
   $g^{(2)}(\partial M_F)$, since no compression move increases the genus of the
 boundary. On the other hand, since $F$ is essential, then $g^{(2)}(\partial M_F)<
    g^{(2)}(\partial M)$, which implies $c(N)<c(M)$.

  Using the inductive assumption, we may conclude that $N$ has a unique root. The
  same is true for $M_F$ and $M_G$, since by Lemma~\ref{lm:decrease} their
  complexities are also smaller than  $c(M)$. Now we have:
  \begin{enumerate}
  \item $M_F$ and $N$ have the same root (since they have a common root
  by Lemma~\ref{commonroot}).
  \item $M_G$ and $N$ have the same root (same reason);
  \item Hence $M_F$ and $M_G$ have the same root, which is a contradiction.
  \end{enumerate}

 \section{Other roots}

 {\sc Roots of cobordisms.}  Recall that a {\em 3-cobordism} is a triple
 $(M,\partial_-M,\partial_+M)$, where $M$ is a compact 3-manifold and
 $\partial_-M$, $\partial_+M$ are unions of connected components of $\partial M$
 such that $\partial_-M\cap \partial_+M=\emptyset$ and $\partial_-M\cup \partial_+M=\partial M$. One can
 define $S$- and $D$-moves on cobordisms just in the same way as
 for manifolds. The $A$-move on cobordisms differs from
 the one for manifolds only in that one boundary circle of $A$
     must lie in $\partial_-M$ while the other   in $\partial_+M$.

\begin{theorem} \label{th:rootsco} For any compact 3-cobordism
$(M,\partial_-M,\partial_+M)$    its
 root   exists and is unique up to homeomorphisms of cobordisms and
removing disjoint 3-spheres and balls.
\end{theorem}

The proof of this theorem is the same as the proof of Theorem~\ref{th:roots}.

{\sc $S$-Roots of manifolds.} We define an $S$-root of $M$ as a
manifold which can be obtained from $M$ by essential $S$-moves and
does not admit any further essential $S$-moves.

\begin{theorem} \label{th:S-roots} For any compact 3-manifold
$ M $,    its
 $S$-root   exists and is unique up to homeomorphism   and
removing disjoint 3-spheres.
\end{theorem}

This theorem is actually equivalent to the     theorem on  the unique
  decomposition into a connected sum of prime factors. Indeed, the $S$-root of
  $M$ coincides with the union of the irreducible prime factors of $M$.

{\sc $(S,D)$-Roots of manifolds.}  An $(S,D)$-root of $M$ is a
manifold which can be obtained from $M$ by essential $S$- and $D$-
moves and does not admit any further essential $S$-moves and
$D$-moves.

\begin{theorem} \label{th:SD-roots} For any compact 3-manifold
$ M $    its
 $(S,D)$-root   exists and is unique up to homeomorphism   and
removing disjoint 3-spheres and balls.
\end{theorem}

For irreducible manifolds this theorem can be deduced from    the  theorem 
 of F. Bonahon~\cite{Bo}
on   characteristic
compression bodies as well as from  \cite{Ma}, where  $D$-roots of irreducible
manifolds     had been  
 considered     under the name {\em cores}.

Our way for proving Theorem~\ref{th:roots} works also for $S$-  and
 $(S,D)$-roots. All   we need is to forget about discs and annuli in the
 first case and about annuli in the second. This  makes the
 proof  significantly shorter.

\end{document}